\documentclass{article}

\def\ra{\rightarrow}

\def\ss{\subseteq}

\def\dbar{\overline{\partial}}

 \def\HollowBox #1#2{{\dimen0=#1 \advance\dimen0 by -#2       
       \dimen1=#1 \advance\dimen1 by #2                       
        \vrule height #1 depth #2 width #2                    
        \vrule height 0pt depth #2 width #1                   
        \llap{\vrule height #1 depth -\dimen0 width \dimen1}%
       \hskip -#2                                             
       \vrule height #1 depth #2 width #2}}                   
 \def\BoxOpTwo{\mathord{\HollowBox{6pt}{.4pt}}\;}             

\def\endpf{\hfill $\BoxOpTwo$}

\font\teneufm=eufm10
\font\seveneufm=eufm7
\font\fiveeufm=eufm5
\newfam\eufmfam
\textfont\eufmfam=\teneufm
\scriptfont\eufmfam=\seveneufm
\scriptscriptfont\eufmfam=\fiveeufm

\newfam\msbfam
\font\tenmsb=msbm10  \textfont\msbfam=\tenmsb
\font\sevenmsb=msbm7  \scriptfont\msbfam=\sevenmsb
\font\fivemsb=msbm5    \scriptscriptfont\msbfam=\fivemsb
\def\Bbb{\fam\msbfam \tenmsb}

\def\CC{{\Bbb C}}

\newtheorem{theorem}{Theorem}
\newtheorem{corollary}[theorem]{Corollary}
\newtheorem{proposition}[theorem]{Proposition}
\newtheorem{lemma}[theorem]{Lemma}
\newtheorem{remark}[theorem]{Remark}

\newtheorem{definition}{Definition}
\newtheorem{example}[definition]{EXAMPLE}

\makeindex

\usepackage{graphicx}

\usepackage{amsmath}

\begin{document}

\begin{center}
\huge \bf
Normed Domains of Holomorphy\footnote{{\bf Key Words:}  domain of holomorphy,
bounded holomorphic functions, normed domain of holomorphy}\footnote{{\bf MR Classification
Numbers:}  32T05, 32E40, 32T25}\end{center}
\vspace*{.12in}

\begin{center}
\large Steven G. Krantz\footnote{Author supported in part
by the National Science Foundation and by the Dean of the Graduate
School at Washington University.}
\end{center}
\vspace*{.15in}

\begin{center}
\today
\end{center}
\vspace*{.2in}

\begin{quotation}
{\bf Abstract:} \sl
We treat the classical concept of domain of holomorphy
in $\CC^n$ when the holomorphic functions considered
are restricted to lie in some Banach space.  Positive
and negative results are presented.  A new view
of the case $n = 1$ is considered.
\end{quotation}
\vspace*{.25in}

\setcounter{section}{-1}

\section{Introduction}

In this paper a {\it domain} $\Omega \ss \CC^n$ is a connected open set.
We let ${\cal O}(\Omega)$ denote the algebra of holomorphic functions
on $\Omega$.

We shall use the following notation:  $D$ denotes
the unit disc in the complex plane.  We let $D^2 = D \times D$ denote the bidisc,
and $D^n = D \times D \times \cdots \times D$ the polydisc in $\CC^n$.  The symbol
$B = B_n $ is the unit ball in $\CC^n$. 

A domain $\Omega \ss \CC^n$ is said to be {\it Runge} if any holomorphic $f$ on $\Omega$
is the limit, uniformly on compact subsets of $\Omega$, of polynomials.

In the classical function theory of several complex variables there are
two fundamental concepts:  domain of holomorphy and pseudoconvex domain.
The Levi problem, which was solved comprehensively in the 1940s and 1950s,
asserts that these two concepts are equivalent:  A domain $\Omega \ss \CC^n$ is
a domain of holomorphy if and only if it is pseudoconvex.  These matters
are discussed in some detail in [KRA1].

Roughly speaking, if $\Omega$ is a domain of holomorphy, then there is
a holomorphic function $f$ on $\Omega$ such that $f$ cannot be analytically
continued to any larger domain.  Generally speaking one cannot say much
about the nature of this $f$---whether it is bounded, or satisfies some other 
growth condition.

In the paper [SIB1], N. Sibony presents the following remarkable example:

\begin{example} \rm
There is a bounded pseudoconvex Runge domain $\Omega \ss \CC^2$, with $\Omega$
a proper subset of the bidisc $D^2 = D \times D$, such that any
{\it bounded} holomorphic function $\varphi$ on $\Omega$ analytically
continues to all of $D^2$.
\end{example}

Of course this result of Sibony can be extended to $\CC^n$ in a variety
of ways.  For one thing, one may take the product of the Sibony domain
with the polydisc in $\CC^{n-2}$ to obtain an example in $\CC^n$.  Alternatively,
one may replace the first (or $z$) variable in the Sibony construction with
a tuple in $\CC^{n-1}$ to obtain a counterexample in $\CC^n$.

It is interesting to note that, in some sense, the Sibony example is generic.  In
fact we have:

\begin{proposition} \sl
The collection of domains $\Omega \ss D^2$, with $\Omega \ne D^2$, such
that any bounded holomorphic function on $\Omega$ analytically continues
to all of $D^2$ (as in the Sibony example above) is uncountable.
\end{proposition}
{\bf Proof:}  We very quickly review the key steps of the Sibony construction.

Let $p_j$ be a countable collection of points in the unit disc $D$, with no
interior accumulation point, so that every boundary point of $D$ is an accumulation
point of the set $\{p_j\}$.  Now define
$$
\varphi(\zeta) = \sum_j \lambda_j \log \left | \frac{\zeta - p_j}{2} \right | \, .
$$
Here $\{\lambda_j\}$ is a summable sequence of positive, real numbers.  Notice that
the function $\varphi$---being the sum of subharmonic functions---is subharmonic.
Define
$$
V_0(\zeta) = \exp (\varphi(\zeta)) \, .
$$
Then $V_0$ has the properties:
\begin{itemize}
\item $V_0$ is subharmonic;
\item $0 < V_0(\zeta) \leq 1$ for all $\zeta \in D$;
\item The function $V_0$ is continuous.
\end{itemize}
The last property holds because the sequence $\{p_j\}$ is discrete and $V_0$ takes
the value 0 only at the $p_j$.

Now define the domain
$$
M(D, V_0) = \{(z, w) \in \CC^2: z \in D, w \in \CC, |w| < \exp(-V_0(z)) \} \, .
$$
Since $V_0$ is positive, we see that this definition makes sense and
that $M(D, V_0)$ is a proper subset of $D^2$.  

The remainder of Sibony's argument shows that any bounded, holomorphic function
on $M(D, V_0)$ analytically continues to a bounded, holomorphic function
on $D^2$.   We shall not repeat it, but refer the reader to [SIB].
  
The key fact in the Sibony construction is that the points $\{p_j\}$ form a discrete
set that accumulates at every boundary point of $D$.  Apart from this property, there
is complete freedom in choosing the $p_j$.   We begin by showing how to construct
two biholomorphically distinct instances of Sibony domains, and then consider at
the end how to produce uncountably many biholomorphically distinct domains.

Define, for $\ell = 1, 2, \dots$,
$$
S_\ell = \{z \in D: |z| = 1 - 2^{- \ell -1}\} \, .
$$
Set 
$$
\{p_j^1\} = \hbox{(the sequence consisting of 4 equally spaced points on}\ S_1 \ , 
$$
$$
                    \hbox{8 equally spaced points on}\ S_2 \ , \
$$
$$
		    \hbox{16 equally spaced points on} \ S_3 \ , \ \hbox{etc.}) 
$$		     
and
$$
\{p_j^2\} = \hbox{(the sequence consisting of 8 equally spaced points on}\ S_1 \ , 
$$
$$
                    \hbox{16 equally spaced points on}\ S_2 \ , \
$$
$$
		    \hbox{32 equally spaced points on} \ S_3 \ , \ \hbox{etc.}) 
$$
Define a domain $\Omega_1$ using the Sibony construction, as above, with
the sequence $\{p_j^1\}$ and define a domain $\Omega_2$ using the Sibony construction 
with the sequence $\{p_j^2\}$.  We claim that $\Omega_1$ and $\Omega_2$ are
biholomorphically inequivalent.

To see this, suppose the contrary. So there is a biholomorphic
mapping $\Phi: \Omega_1 \rightarrow \Omega_2$. By the usual
classical arguments (see the proof of Proposition 11.1.2 in
[KRA1]), we see that $\Phi$ must commute with rotations in the
$w$ variable.	It follows that any disc in $\Omega_1$ of the form
$$
{\bf d}_z = \{(z, w): |w| < \exp(-V_0(z))\} \, ,
$$
for $z \in D$ fixed, must be mapped to a similar disc in $\Omega_2$.

Further observe that each of the discs ${\bf d}_{p_j} \ss \Omega_1$ is a totally
geodesic submanifold in the Kobayashi metric.  This assertion follows
immediately from the existence of the maps
$$
D \stackrel{i}{\longrightarrow} \Omega_1 \stackrel{\pi}{\longrightarrow} D \, ,
$$
where $i$ is the injection
$$
i(w) = (p_j,w)
$$
and $\pi$ is the projection
$$
\pi(p_j,w) = w \, .
$$
Observe that $\pi \circ i = \hbox{id}$.    Similar reasoning shows that the disc ${\bf d}^* = \{(z,0)\} \ss \Omega_1$
is totally geodisc.  Of course similar remarks apply to the corresponding discs in $\Omega_2$.

Now it is essential to notice that, in $\Omega_1$, the vertical discs of the form ${\bf d}_z$ for $z$ not
one of the $p_j$ are {\it not} totally geodesic.  This follows because, at the point $(z,0)$, the Kobayashi extremal
disc in the vertical direction $\langle 0, 1 \rangle$ will not be the rigid disc $\zeta \mapsto (0, \zeta)$ but
rather a disc that curves into one of the spikes above a nearby $p_j$.

A final observation that we need to make is this.    The vertical totally geodesic discs will be mapped to
each other by the biholomorphic mapping $\Phi$.  But more is true.  Because ${\bf d}^*$ is totally geodesic,
in fact the totally geodesic discs located at the points $p_j$ that lie on $S_1$ in $\Omega_1$ must be mapped
to the totally geodesic discs located at the points $p_j$ that lie on $S_1$ in $\Omega_2$ (because both
circles consist of points that have the same Kobayashi distance from the origin).  This is impossible because
$S_1$ in $\Omega_1$ has four such points while $S_1$ in $\Omega_2$ has eight such points.   That is the
required contradiction.

Now it is clear how to construct uncountably many inequivalent
domains of Sibony type. For if ${\cal A} = \{\alpha_j\}$ is
any exponentially increasing sequence of positive integers
then we may associate to it a domain $\widetilde{\Omega}$ of
Sibony-type with $\alpha_1$ points on $S_1$, $\alpha_2$ points
on $S_2$, and so forth. The preceding argument shows that
different choices of ${\cal A}$ result in biholomorphically
inequivalent domains. And there are clearly uncountably many
such sequences. That completes the argument. 
\endpf 
\smallskip \\
       
The Sibony result has an interesting and important interpretation in terms
of the corona problem.  We have:

\begin{proposition} \sl
Let $\Omega \ss \CC^n$ be a bounded domain.  Suppose
that $X$ is a Banach space of holomorphic functions
on $\Omega$ that contains $H^\infty(\Omega)$.  Let $\Omega'$ be a strictly larger domain that
contains $\Omega$.  Assume that any element of $X$ analytically
continues to a holomorphic function on $\Omega'$ (we often assume that
the extended function satisfies a similar norm estimate to that specified
by the norm on $X$, but that is not necessary and we do not impose
that condition at this time).  Then the corona problem cannot
be solved in the space $X$.  That is to say, if $f_1, f_2, \dots, f_k$ are
holomorphic functions in $X$ with no common zero, then there do {\it not}
exists elements $g_1, g_2, \dots, g_k \in X$ such that
$$
f_1 g_1 + f_2 g_2 + \cdots f_k g_k \equiv 1
$$
on $\Omega$.
\end{proposition}
{\bf Proof:}   Assume to the contrary that such $g_1, g_2, \dots, g_k$ exist.
Of course each $g_j$ analytically continues to $\Omega'$.  

Let $P = (p_1, p_2, \dots, p_n)$ be a point of $\Omega' \setminus \Omega$.   Set
$f_j(z) = z_j - p_j$.  Then the $f_j$ have no common zero in $\Omega$.  So,
by hypothesis, the $g_j$ exist.  And these functions extend analytically
to $\Omega'$.  But then
$$
f_1 g_1 + f_1 g_2 + \cdots f_n g_n \equiv 1
$$
on $\Omega'$.  Since the $f_j$ all vanish at $P$, we see that, at $P$, 
the lefthand side of this last equation vanishes.  That is clearly a contradiction.
Hence the $g_j$ do not exist.
\endpf 
\smallskip \\

Of course this last proposition means in particular that the point evaluations
on $\Omega$ are not weak-$*$ dense in the maximal ideal space
of $H^\infty(\Omega)$.  See [KRA3] for more on these matters.

By contrast to Sibony's result, David Catlin [CAT1] has shown that any smoothly bounded, pseudoconvex
domain in $\CC^n$ supports a bounded holomorphic function that cannot
be analytically continued to any larger domain.   In fact he has proved
something sharper:

\begin{theorem} \sl
Let $\Omega \ss \CC^n$ be a smoothly bounded pseudoconvex domain.  Then
there is a function in $C^\infty(\overline{\Omega})$, holomorphic on
the interior, which cannot be analytically continued to any larger domain.
\end{theorem}

Hakim and Sibony [HAS] have proved something even more decisive:

\begin{theorem} \sl
Let $\Omega \ss \CC^n$ be a smoothly bounded pseudoconvex domain.  Then
the maximal ideal space (or spectrum) of the algebra $C^\infty(\overline{\Omega}) \cap {\cal O}(\Omega)$
is in fact $\overline{\Omega}$.
\end{theorem}

It should be stressed that the proofs of the last two results use an algebraic formalism of 
H\"{o}rmander [HOR] which entails the loss of some derivatives,
so it is essential to be working with functions that are $C^\infty$ on $\overline{\Omega}$.  Attempts
to adapt the arguments to other function spaces are doomed to failure.

Peter Pflug [PFL] has shown that the situation for $L^2$ holomorphic functions is very
neat and elegant:

\begin{theorem} \sl
Let $\Omega \ss \CC^n$ be {\it any} pseudoconvex domain.  Then there
is an $L^2$ holomorphic function on $\Omega$ that cannot be analytically
continued to a larger domain.
\end{theorem}

It is natural to ask for a characterization of those domains $\Omega$ which
are domains of holomorphy in the traditional sense but {\it not} domains of
holomorphy for bounded holomorphic functions.  One would also like to know
whether there are analogous results for $L^p$ holomorphic functions, $1 \leq p < \infty$.

The purpose of the present paper is to consider these matters.  While
we cannot provide a full answer to the questions just posed, we can certainly
give some useful partial results, and point in some new directions.

We mention in passing that the paper [DAR] contains some results that bear
on the questions posed here.  The arguments presented in [DAR] appear to
be incomplete.

\section{Some Notation}

Let us say that a domain $\Omega \ss \CC^n$ is of type $HL^p$, $1 \leq p \leq \infty$,
if there is a holomorphic function $f$ on $\Omega$, $f \in L^p(\Omega)$, which
cannot be analytically continued to any larger domain.  We instead say
that $\Omega$ is of type $EL^p$ if there is a strictly larger domain $\widehat{\Omega}$ so that
every holomorphic $L^p$ function on $\Omega$ analytically continues to $\widehat{\Omega}$.
Obviously $HL^p$ and $EL^p$ are disjoint.  Every domain is either of type $HL^p$ or $EL^p$.

We are interested in giving an extrinsic description of those domains which are of type
$HL^p$ and those which are of type $EL^p$.

\section{The Situation in the Complex Plane}

Matters in one complex variable are fairly well understood.

First of all, we should note the example of $$ \Omega = D
\setminus \{0\} \, . $$ Of course, by the Riemann removable
singularities theorem, any bounded holomorphic function on
$\Omega$ analytically continues to all of $D$. So $\Omega$ is
{\it not} a domain of type $HL^\infty$. It {\it is} a domain
of type $EL^\infty$.

In fact it may be noted (for the domain $\Omega$ in the last
paragraph) that, if $p \geq 2$, then any holomorphic function
that is $L^p(\Omega)$ will analytically continue to all of
$D$ (see [HAP]). So this $\Omega$ is a domain of type $EL^p$.  By contrast,
if $p < 2$ then the function $f(\zeta) = 1/\zeta$ is holomorphic on $\Omega$
and in $L^p(\Omega)$.  But of course this $f$ {\it does not} analytically
continue to the full disc $D$.  So, for $p < 2$, the domain is of
type $HL^p$.  

The treatment in [HAP] of the matter just discussed is rather abstract, and
it is worthwhile to have a traditional function-theoretic treatment of
these matters.  We provide one now.  We thank Richard Rochberg for a helpful
conversation about this topic.  So let $f$ be holomorphic on $D \setminus \{0\}$
and assume that $f \in L^2(D)$ (the case $p > 2$ follows immediately from this one).
We write $f(\zeta) = \sum_{j=-\infty}^\infty a_j \zeta^j$.  For $0 < a < b < 1$ and $k$ a negative
integer, consider the expression
$$
A = \int_a^b r \int_0^{2\pi} f(\zeta) e^{-i k \theta} \, d\theta \, dr \, ,
$$
where it is understood that $\zeta = r e^{i\theta}$.

On the one hand,
\begin{eqnarray*}
|A| & \leq & \int_{a \leq |\zeta| \leq b} |f(\zeta)| \, dA(\zeta) \\
    & \leq & \|f \|_{L^2} \cdot \left | \{\zeta: a \leq |\zeta| \leq b\} \right |^{1/2} \\
    & = & \|f\|_{L^2} \cdot \left [ \pi (b^2 - a^2) \right ]^{1/2} \, .
\end{eqnarray*}
On the other hand,
\begin{eqnarray*}
|A| & = & \left | \int_a^b r a_k r^k \, dr \right |  \\
    & = & \left | a_k \cdot \frac{r^{k+2}}{k+2} \biggr ]_a^b \right | \\
    & = & \left | a_k \left ( \frac{b^{k+2}}{k+2} - \frac{a^{k+2}}{k+2} \right ) \right | \, .
\end{eqnarray*}
If $k < -2$ and $a = b/2$, this gives a contradiction as $b \ra 0^+$.  Of course the
cases $k = -2$ and $k = -1$ can be handled separately because $\zeta^{-2}$ and $\zeta^{-1}$ are
certainly {\it not} in $L^2$.
\endpf 
\smallskip \\

\begin{remark} \rm
We note that the proof goes through for $p < 2$ up until the very end.  One must
note that $\zeta^{-1}$ in fact {\it does} lie in $L^p$ for $p < 2$.  So there
is no removable singularities theorem for this range of $p$.
\end{remark}

\begin{remark} \rm
A standard result coming from potential theory is that, if $\Omega \ss \CC$, $P \in \Omega$, and $f$ is holomorphic
on $\Omega \setminus \{P\}$, then $|f(z)| = o(\log[1/|z - P|])$ (where we are using Landau's notation) implies that 
$f$ continues analytically to all of $\Omega$.  The philosophy here is that a function $f$ satisfying this
growth hypothesis has a singularity at $P$ that is milder than the singularity of the Green's function.
This point of view is particularly useful in the study of removable singularities for harmonic functions.
This result is not of any particular interest for us because it is not formulated in the language
of Lebesgue spaces.  In any event, it is weaker than the result presented above for $L^2$ because the
logarithm function is certainly square integrable.   It is a pleasure to thank Al Baernstein and David Minda for helpful
remarks about these ideas.
\end{remark}

The enemy in the results discussed at the beginning of this section is that $\Omega$ is not equal to the interior of its closure.
In fact we have:

\begin{proposition} \sl
Suppose that the bounded domain $\Omega \ss \CC$ is the interior of its closure.  Then $\Omega$ is
a domain of type $HL^p$ for $1 \leq p \leq \infty$.
\end{proposition}
{\bf Proof:}  The proof that we now present is an adaptation and simplification 
of an argument from [BER].

Let $\{p_j\}$ be a countable, dense subset of ${}^c \overline{\Omega}$.  For each
$j$, the function $\varphi_j(\zeta) = 1/(\zeta - p_j)$ is holomorphic and
bounded on $\Omega$ and does not analytically continue past $p_j$.  

Now, for each $j$, let ${\bf d}_j$ be an open disc centered at $p_j$ which has
nontrivial intersection with $\Omega$.  Consider the linear mapping
$$
I_j: {\cal O}(\Omega \cup {\bf d}_j) \cap L^p(\Omega \cup {\bf d}_j) \longrightarrow {\cal O}(\Omega) \cap L^p(\Omega)
$$
given by restriction.  Of course each of the indicated spaces is equipped with the $L^p$ norm, and is therefore
a Banach space.  We note that the example above of $\varphi_j$ shows that $I_j$ is not surjective.
As a result, the open mapping principle tells us that the image ${\cal M}_j$ of $I_j$ is of first category
in ${\cal O}(\Omega) \cap L^p(\Omega)$.
Therefore, by the Baire category theorem,
$$
{\cal M} \equiv \bigcup_j {\cal M}_j
$$
is of first category in ${\cal O}(\Omega) \cap L^p(\Omega)$.  But this just says that
the set of $L^p$ holomorphic functions on $\Omega$ that can be analytically continued to some $p_j$ is
of first category.  Therefore the set of $L^p$ holomorphic functions that {\it cannot} be
analytically continued across the boundary is dense in ${\cal O}(\Omega) \cap L^p(\Omega)$.  That completes the proof.
\endpf
\smallskip \\

The key point of the proof just presented is that, for each point not in the closure
of the given domain, there is a function holomorphic on the domain (and in the
given function space) that does not analytically continue past the point.  Such
functions are trivial to construct in one complex variable.  Not so in higher
dimensions.

We note in passing that, when $\Omega \ss \CC$ is the unit disc $D$ then it is easy
to construct a bounded holomorphic function that does not analytically continue
to a larger domain.  For let $\{p_j\}$ be a discrete set in $D$ that accumulates
at every boundary point and so that
$$
\sum_j 1 - |p_j| < \infty \, .
$$
For example, take
\begin{eqnarray*}
 && p_1, p_2, p_3, p_4 \ \ \hbox{to be equally spaced points at distance 1/4 from $\partial D$} \\
 && p_5, p_6, \dots, p_{12} \ \ \hbox{to be equally spaced points at distance 1/8 from $\partial D$}  \\
 && p_{13}, p_6, \dots, p_{28} \ \ \hbox{to be equally spaced points at distance 1/16 from $\partial D$}  \\
\end{eqnarray*}
and so forth.  Then the Blaschke product with zeros at the $p_j$ will do the job.   If $\Omega$ is
a simply connected domain having a Jordan curve as its boundary, then conformal mapping together with
Carath\'{e}odory's theorem about continuous boundary extension will give a bounded, holomorphic, non-continuable
function on this $\Omega$.

We close this section by noting that, if $\Omega$ is a domain of holomorphy in $\CC^n$ and if 
$V = \{z \in \Omega: f(z) = 0\}$ for some $f$ holomorphic $f$ on $\Omega$ (we call
$V$ a {\it variety}) then $\Omega' = \Omega \setminus V$ is also a domain of holomorphy (for if $\varphi$
is a holomorphic function on $\Omega$ that does not analytically continue to a larger domain
then $\varphi/f$ is a holomorphic function on $\Omega'$ that does not analytically continue
to any larger domain.   And it is easy to see that $\Omega'$ is an $EL^\infty$ domain.
See [GUR, p.\ 19] for the details.  

\section{Complications in Dimension \boldmath $n$}

As we have indicated, the example of Sibony exhibits a domain
which is {\it not} of type $HL^\infty$ (instead it is
of type $EL^\infty$).  The theorem of Catlin shows 
that all smoothly bounded, pseudoconvex domains are
of type $HL^\infty$.

It of course makes sense to focus this discussion on pseudoconvex
domains.  If a domain $\Omega$ is {\it not} pseudoconvex then there
will perforce be a larger domain $\Omega'$ to which all holomorphic
functions (regardless of growth) on $\Omega$ analytically continue.
So this situation is not interesting.

Thus we see that the domains of interest for us will be pseudoconvex
domains that do {\it not} have smooth boundary.  Our first result
is as follows:

\begin{proposition} \sl
Let $D_1, D_2, \dots, D_n$ be bounded domains in $\CC$, each of which
is equal to the interior of its closure.  Define
$$
\Omega = D_1 \times D_2 \times \cdots \times D_n \, .
$$
Then $\Omega$ is a domain of type $HL^p$ for any $1 \leq p \leq \infty$.
\end{proposition}
{\bf Proof:}  Fix $p$ as indicated.  Then, by Proposition 8, there
is a holomorphic function $\psi_j$ on $D_j$, for $1 \leq j \leq n$, such
that $\psi_j$ is holomorphic and $L^p$ on $D_j$ and does not
analytically continue to any larger domain.

But then
$$
\psi(z_1, \dots, z_n) = \psi_1(z_1) \cdot \psi_2(z_2) \cdot \cdots \cdot \psi_n(z_n)
$$
is holomorphic and $L^p$ on $\Omega$ and does not analytically continue
to any larger domain.
\endpf 
\smallskip \\

\begin{proposition} \sl
Let $\Omega \ss \CC^n$ be bounded and convex.  Let $1 \leq p \leq \infty$.
Then $\Omega$ is a domain of type $HL^p$.
\end{proposition}
{\bf Proof:}  Just imitate the proof of Proposition 8.  The main point
to note is that if $q \not \in \overline{\Omega}$ and $\nu$ is a unit normal
vector from $\partial \Omega$ out through $q$, then the function
$$
\psi(z) = \frac{1}{(z - q) \cdot \nu}
$$
is holomorphic and bounded on $\Omega$ and is singular at $q$.  So the
rest of the proof goes through as before.
\endpf 
\smallskip \\

In fact more is true:

\begin{proposition} \sl Let $\Omega \ss \CC^n$ be bounded and
strongly pseudoconvex with $C^2$ boundary. Let $1 \leq p \leq
\infty$. Then $\Omega$ is a domain of type $HL^p$.
\end{proposition}
{\bf Proof:}  Of course we again endeavor to apply the argument of the proof
of Proposition 8.   It is enough to restrict attention to points
$q$ in ${}^c \Omega$ which are sufficiently close to $\partial \Omega$.
If $q$ is such a point then there is a larger strongly pseudoconvex domain $\Omega'$
with $C^2$ boundary such that $\overline{\Omega} \ss \Omega'$ and $q \in \partial \Omega'$.
Now let $L_q(z)$ be the Levi polynomial (see [KRA1]) for $\Omega'$ at $q$.  Then there
is a neighborhood $U$ of $q$ so that 
$$
\{z \in U: L_q(z) = 0\} \cap \overline{\Omega'} = \{q\} \, .
$$

Thus $f_q(z) = 1/L_q(z)$ is holomorphic on $\Omega' \cap U$ and singular at $q$.
Let $\varphi$ be a $C_c^\infty$ function that is compactly supported in $U$ and
is identically equal to 1 in a small neighborhood of $q$.  We wish to choose
a bounded function $h$ so that
$$
g(z) = \frac{\varphi(z)}{L_q(z)} + h
$$
is holomorphic on $\Omega'$.  This entails solving the $\overline{\partial}$-problem
$$
\dbar h = - \frac{\dbar \varphi(z)}{L_q(z)}  \, .
$$
Of course the data on the righthand side of this equation is $\dbar$-closed with bounded
coefficients.  By work in [KRA2] or [HEL] we see that a bounded solution $h$ exists.

This gives us a functions $g$ that is {\bf (i)} holomorphic on $\Omega'$ and {\bf (ii)}
singular at $q$.   This is just what we need, for points $q$ in ${}^c \Omega$ that are 
close to $\partial \Omega$, in order to imitate the proof of Proposition 8.   That completes
our argument.	  See also [BER, Theorem 3.6] for a similar result with a somewhat different
proof in the case $p = \infty$.
\endpf
\smallskip \\

For finite type domains we can prove the following result.    Let $\Omega \ss \CC^2$ be
given by $\Omega = \{z \in \CC^2: \rho(z) < 0\}$.  Recall that a point
$q$ in the boundary of a domain $\Omega$ is said to be of {\it finite
geometric type} $m$ in the sense of Kohn if there is a nonsingular, one-dimensional analytic
variety $\varphi: D \ra \CC^2$ with $\varphi(0) = q$ and
$$
|\rho(\varphi(\zeta))| \leq C \cdot |\varphi(\zeta) - q|^m \, ,
$$
and so that there is no other nonsingular, one-dimensional analytic
variety satisfying a similar inequality with $m$ replaced by $m + 1$.   These
ideas are discussed in detail in [KRA1, Ch.\ 10].

It is known that the geometric definition of finite type given in the last paragraph
is equivalent to a more analytic one in terms of commutators of vector fields.  Namely,
let 
$$
L = \frac{\partial \rho}{\partial z_1} \frac{\partial}{\partial z_2} -  \frac{\partial \rho}{\partial z_2} \frac{\partial}{\partial z_1}
$$
be a complex tangential vector field to $\partial \Omega$ and $\overline{L}$ its conjugate.	 A {\it first order commutator} is a Lie bracket
of the form $[L, \overline{L}]$.  A {\it second order commutator} is a Lie bracket of the form
$[L, M]$ or $[\overline{L}, M]$, where $M$ is a first order commutator.  And so forth.  We say
that a point $q \in \partial \Omega$ is of {\it analytic type} $m$ if all the commutators
${\cal L}$ up to and including order $m - 1$ have the property that 
$$
{\cal L}(\rho)[q] = 0
$$
but there is a commutator ${\cal L}'$ of order $m$ such that
$$
{\cal L}' (\rho)[q] \ne 0 \, .
$$ 

It is a result of Kohn [KOH] and Bloom/Graham [BLG] that, when $\Omega \ss \CC^2$, a point $q \in \partial \Omega$ is of geometric
finite type if and only if it is of analytic finite type.   Details of these matters may be found
in [KRA1].

Now it is easy to see that the notion of analytic finite type varies semi-continuously with smooth variation
of $\rho$.   In particular, if each point of $\partial \Omega$ is of some finite type, then the type
of the point will vary semi-continously.  So there is an upper bound $M$ for all types of points
in $\partial \Omega$.  In this circumstance we say that {$\Omega$ is a domain of finite type} (at most) $M$.
		     
As a result of these considerations, we have the following lemma:

\begin{lemma} \sl
Let $\Omega = \{z \in \CC^2: \rho(z) < 0\}$ be a domain of finite type $M$.  Then there are domains $\Omega'$ of finite type
so that $\Omega' \supset \overline{\Omega}$.  In particular, if $\psi$ is a smooth, negative
function with $\|\psi\|_{C^{M+1}}$ sufficiently small and $\rho' = \rho + \psi$ then
$\Omega' \equiv \{z \in \CC^2: \rho'(z) < 0\}$ will contain $\overline{\Omega}$ and be of finite type.
\end{lemma}

Now we have

\begin{proposition} \sl Let $\Omega \ss \CC^2$ be smoothly bounded and
of finite type $m$. Let $1 \leq p \leq
\infty$. Then $\Omega$ is a domain of type $HL^p$.
\end{proposition}
{\bf Proof:}  The argument is similar to that for the last few propositions.  If $q \not \in \overline{\Omega}$
and is sufficiently close to $\partial \Omega$, then we may use the last lemma and the discussion preceding
that to construct a finite type domain $\Omega' \supset \overline{\Omega}$ and with $q \in \partial \Omega'$.
Now the theorem of Bedford and Forn\ae ss [BEF] gives us a peaking function $f_q$ for the point $q$ on
the domain $\Omega'$.   That is to say, 
\begin{enumerate}
\item[{\bf (i)}]  $f_q$ is continuous on $\overline{\Omega'}$;
\item[{\bf (ii)}]  $f_q$ is holomorphic on $\Omega'$;
\item[{\bf (iii)}]  $|f_q(z)| \leq 1$ for all $z \in \overline{\Omega'}$;
\item[{\bf (iv)}]  $f_q(q) = 1$;
\item[{\bf (v)}]  $|f_q(z)| < 1$ for all $z \in \overline{\Omega'} \setminus \{q\}$.
\end{enumerate}
Then the function $g_q(z) = 1/[1 - f_q(z)]$ is holomorphic on $\Omega$ and singular
at $q$.  

The rest of the argument is completed as in the proof of the last Proposition.
\endpf 
\smallskip \\
	     
We note that the Kohn-Nirenberg domain [KON] shows that, even on a finite type
domain in $\CC^2$, we cannot hope for a holomorphic separating function
like $L_q$ in the strongly pseudoconvex case.  But the peak function
of Bedford-Forn\ae ss suffices for our purposes.

\begin{proposition} \sl
Let $\Omega \ss \CC^n$ be a bounded analytic polyhedron.  Certainly $\Omega$
is then a domain of holomorphy.	 We have that $\Omega$ is
a domain of type $HL^p$ for $1 \leq p \leq \infty$.
\end{proposition}
{\bf Proof:}  We know by the standard definition (see [KRA1]) that
$$
\Omega = \{z \in \CC^n :|f_1(z)| < 1, |f_2(z)| < 1, \dots, |f_k(z)| < 1\}
$$
for some holomorphic functions $f_j$.  Now if $q \not \in \overline{\Omega}$, 
then there is some complex constant $\lambda$ with $|\lambda| > 1$ and some $j$
so that $f_j(q) = \lambda$.  That being the case, the function
$$
\psi(z) = \frac{1}{\lambda - f_j(z)}
$$
is a function that is bounded and holomorphic on $\Omega$ but singular at $q$.
Now the proof can be completed as in the previous Propositions.
\endpf 
\smallskip \\

\begin{proposition} \sl 
Let $\Omega \ss \CC^n$ be a complete circular domain.  Assume that $\Omega$ is pseudoconvex.
Then $\Omega$
is a domain of type $HL^p$, $1 \leq p \leq \infty$.
\end{proposition}
{\bf Proof:}  	Let $q$ be a point that does not lie in $\overline{\Omega}$.
Let $q^*$ be the nearest point to $q$ in the boundary of $\Omega$, and
let $\nu$ be the unit outward normal vector at $q^*$.  Set
$$
f_q(z) = (z - q) \cdot \nu \, .
$$
Then $f_q$ is holomorphic, and we claim that the zero
set ${\cal Z}_q$ of $f_q$ does not intersect $\overline{\Omega}$.
Suppose to the contrary that it does.

Let $x$ be a point that lies in both $\overline{\Omega}$ and
in ${\cal Z}_q$.   Of course any point that can be obtained
by rotating the coordinates of $x$ will also lie in $\overline{\Omega}$.
One such choice of rotations will give a point that lies on the 
ray from the origin out to $q$.  But that rotated point will be further
from the origin than $q$ itself (by the Pythagorean theorem).  Since it lies in $\overline{\Omega}$ then
so does $q$ (because the domain is complete circular).  That is a contradiction.  
Therefore $x$ does not exist and $\overline{\Omega}$ and the
zero set of $f_q$ are disjoint.

As a result, the function $g_q \equiv 1/f_q$ is holomorphic and bounded
on $\Omega$ and singular at $q$.  The  proof may now be completed as
in the preceding propositions.
\endpf 
\smallskip \\

The next result points in the general direction that any reasonable pseudoconvex domain
will be of type $HL^p$ for $1 \leq p \leq \infty$.

\begin{proposition} \sl
Let $\Omega$ be a bounded, pseudoconvex domain with a Stein neighborhood basis.\footnote{Here
a {\it Stein neighborhood basis} for $\Omega$ is a decreasing collection of
pseudoconvex domains $\Omega_j$ such that $\cap_j \Omega_j = \overline{\Omega}$.  See [CHS]
for further details in this matter.}
Then $\Omega$ is a domain of type $HL^p$ for $1 \leq p \leq \infty$.
\end{proposition}

\begin{remark}  \rm
Of course a domain with Stein neighborhood basis can have rough boundary.  So
this proposition says something new and with content.
\end{remark}

\noindent {\bf Proof of the Proposition:}  Let $\epsilon > 0$.  By definition
of Stein neighborhood basis, there is a pseudoconvex domain $\widetilde{\Omega}$ so that
$\widetilde{\Omega} \supseteq \overline{\Omega}$.  Therefore (see [KRA1, Ch.\ 3])
there is a smoothly bounded, strongly pseudoconvex domain $\widetilde{\widetilde{\Omega}}$
so that $\widetilde{\Omega} \supseteq \widetilde{\widetilde{\Omega}} \supset \overline{\Omega}$.
Let $P \in \partial \widetilde{\widetilde{\Omega}}$.  Then we may imitate the construction
in the proof of Proposition 8 to find a function that is holomorphic and bounded on
$\Omega$, extends past the boundary of $\Omega$, but is singular at $P$.   Now
the rest of the argument---elementary functional analysis---is just as in the
proof of Proposition 8.
\endpf 
\smallskip \\

The interest of Propositions 14, 15, and 16 is that the domains
constructed there have only Lipschitz boundary. We know for
certain (thanks to Catlin and Hakim/Sibony) that pseudoconvex domains with
smooth boundary are of type $HL^\infty$. And there are domains
with rough boundary, such as the Sibony domain, that are of
type $EL^\infty$. So the last two Propositions give examples
of domains with rough boundary which are of type $HL^\infty$.

\section{Other Properties of \boldmath $HL^p$ and $EL^p$ Domains}

in [BER] an example is given which shows that the increasing
union of $HL^\infty$ domains need not be $HL^\infty$. Indeed,
it is well known (see [BERS]) that {\it any} domain of
holomorphy is the increasing union of analytic polyhedra (ref.\ Proposition 10). 
Of course an analytic polyhedron is $HL^p$ for $1 \leq p \leq
\infty$, but the Sibony domain (which is certainly the union of
analytic polyhedra) described above is pseudoconvex
and not $HL^\infty$. Berg in addition shows (see his Theorem
1.15) that the decreasing intersection of $HL^\infty$ domains
is $HL^\infty$.

Now we describe some other related examples.  Again see [BER] for cognate ideas.

\begin{example} \rm
The decreasing intersection of $HL^p$ domains is $HL^p$, $1 \leq p \leq \infty$.
To see this, let
$$
\Omega_1 \supseteq \Omega_2 \supseteq \cdots \supseteq 
$$
and $\Omega_0 = \cap_j \Omega_j$.  Assume that each $\Omega_j$ is $HL^p$.
If $\Omega_0$ is not $HL^p$ then it is of course $EL^p$.  So there is
a strict neighborhood $\Omega'_0$ of $\Omega_0$ so that every holomorphic
function in $L^p$ of $\Omega_0$ analytically continues to $\Omega'_0$.
But then certainly any holomorphic function on $\Omega_j$ which is
in $L^p$ will analytically continue to $\Omega'_0$.  When $j$ is large,
this will contradict the fact that $\Omega_j$ is $HL^p$.
\end{example}

\begin{example} \rm
There is a decreasing sequence $\Omega_1 \supseteq \Omega_2 \supseteq \cdots$
of $EL^\infty$ domains such that the intersection domain $\Omega_0 \equiv \cap_j \Omega_j$ 
is not $EL^\infty$.

To see this, we follow the construction of [SIB1, p.\ 206].  Let $\{a_j\}$ be a sequence in the 
unit disc $D$ with no interior accumulation point and such that every boundary point of $D$
is the nontangential limit of some subsequence.  Let $\lambda_j$ be a summable sequence
of positive real numbers.  Define, for $\epsilon > 0$ and $z \in D$,
$$
\varphi^\epsilon(z) = \sum_j \epsilon \lambda_j \log \left | \frac{z - a_j}{2} \right | \, .
$$
Then certainly $\varphi^\epsilon$ is subharmonic and negative on $D$.  Further note
that the functions $\varphi^\epsilon$ increase pointwise to the identically 0 function 
as $\epsilon \rightarrow 0^+$.  Now set
$$
V_0^\epsilon(z) = \exp(\varphi^\epsilon(z)) \, .
$$
Then $V_0^\epsilon$ is also subharmonic, $0 \leq V_0^\epsilon < 1$.  The function
takes the value 0 only at the points $a_j$.   

Finally define the domains
$$
M^\epsilon(D, V_0^\epsilon) = \{(z,w) \in \CC^2: z \in D, w \in \CC, |w| < \exp(- V_0^\epsilon(z))\} \, .
$$
Each $M^\epsilon(D, V_0^\epsilon)$ is pseudoconvex.   And the argument of Sibony shows
that it is a domain of type $EL^\infty$.  But notice that the function
$\exp(-V_0^\epsilon(z))$ decreases pointwise to the function that is identically
equal to $1/e$ as $\epsilon \rightarrow 0^+$.  Hence the domains
$M^\epsilon(D, V_0^\epsilon)$ decrease to the bidisc
$\{(z,w): z \in D, |w| < 1/e\}$.  And the latter is a domain
of type $HL^\infty$.

So we have produced a decreasing sequence of $EL^\infty$ domains whose intersection
is $HL^\infty$.
\end{example}

We now give a separate proof, which has independent interest,
of the contrapositive of Proposition 16.

\begin{proposition} \sl
If $\Omega$ is a bounded domain of type $EL^p$, $1 \leq p \leq \infty$ then
$\Omega$ does not have a Stein neighborhood basis.
\end{proposition}
{\bf Proof:}  Suppose that every holomorphic $L^p$ function
on $\Omega$ analytically continues to a larger domain $\widehat{\Omega}$.
Seeking a contradiction, we assume that $\Omega$ has a Stein neighborhood
basis.  Choose a pseudoconvex domain $U \supseteq \overline{\Omega}$ so
that $\widehat{\Omega} \setminus U$ is nonempty.

Now there is some holomorphic function $g$ on $U$ that does not analytically
continue to any larger open domain.  Therefore the restriction of $g$ to
$\Omega$ is a holomorphic $L^p$ function $\widetilde{g}$ on $\Omega$ that analytically
continues to $U$ but no further.  This contradicts the fact that $\widetilde{g}$ must
analytically continue to $\widehat{\Omega}$.  We conclude that $\Omega$ cannot
have a Stein neighborhood basis.
\endpf
\smallskip 

We close with the following useful property of $EL^\infty$ domains:

\begin{proposition} \sl
Let $\Omega$ be a bounded, $EL^\infty$ domain in $\CC^n$, so that
any bounded, holomorphic function $f$ on $\Omega$ analytically continues
to some bounded, holomorphic function $\widehat{f}$ on some $\widehat{\Omega}$.  Let $f$
be a bounded, holomorphic function on $\Omega$ so that $|f|$ is bounded from
$0$ by some $\eta > 0$.  Then $\widehat{f}$ will be nonvanishing.
\end{proposition}
{\bf Proof:}  Of course $g = 1/f$ makes sense on $\Omega$ and is holomorphic and
bounded, so it analytically continues to some bounded, holomorphic function
$\widehat{g}$ on $\widehat{\Omega}$.  But of course $1 \equiv f \cdot g$ analytically
continues to the identically 1 function on $\widehat{\Omega}$.  So we see
that $\widehat{f} \cdot \widehat{g} \equiv 1$ on $\widehat{\Omega}$.  We conclude
then that $\widehat{f}$ cannot vanish.
\endpf 
\smallskip \\

\section{Relationship with the $\dbar$-problem}

In the paper [SIB2], N. Sibony exhibits a smoothly bounded, pseudoconvex
domain on which the equation
$$
\dbar u = f \,,
$$
for $f$ a $\dbar$-closed $(0,1)$ form with bounded coefficents, has
no bounded solution $u$.  This is important information for function
theory, and also for the theory of partial differential equations.

It is natural to speculate that there is some relation between those
domains on which the $\dbar$-equation satisfies uniform estimates
and those domains which are of type $HL^\infty$.  In that vein,
we offer the following result:

\begin{proposition} \sl
Let $\Omega \ss \CC^n$ be a bounded domain which is of finite type $m$ and so
that the $\dbar$-equation $\dbar u = f$ satisfies uniform estimates
on $\Omega$.  That is to say, there is a universal constant $C > 0$ so that, 
given a $\dbar$-closed $(0,1)$ form $f$ with
bounded coefficients there is a solution $u$ to the equation $\dbar u = f$ with
$$
\|u\|_{L^\infty} \leq C \cdot \|f\|_{L^\infty} \, .
$$
Then $\Omega$ is a domain of type $HL^\infty$.
\end{proposition}

\begin{remark} \rm
It is important to notice in this last Proposition that the domain $\Omega$ need not have
$C^\infty$ boundary.  For type 2, it suffices for the boundary to be $C^2$.  For
type $m \geq 2$, it suffices for the boundary to be $C^m$.

It is known that strongly pseudoconvex domains [KRA2], finite type domains
in $\CC^2$ [FOK], and the polydisc [HEN] all satisfy uniform estimates
for the $\dbar$-problem.

\end{remark}
\vspace*{.12in}

\noindent {\bf Proof of the Proposition:} It is known (see, for example
[CAT2], [CAT3] or [DAN1]), [DAN2]) that the strongly pseudoconvex points in
$\partial \Omega$ form an open, dense set.   Let $q \in \partial \Omega$ be such
a point and let $\epsilon > 0$.  Let $\nu$ be the unit outward normal
vector to $\partial \Omega$ at $q$ and set $q' = q + \epsilon \nu$.  If $\epsilon$ 
is small then there is a ``bumped domain'' $\Omega'$ with these properties:
\begin{itemize}
\item There is a small neighborhood $U$ of $q$ so that $U \cap \partial \Omega$ 
consists only of strongly pseudoconvex points.	
\item $\partial \Omega \setminus U = \partial \Omega' \setminus U$.
\item $\partial \Omega' \cap U$ is strongly pseudoconvex and lies outside $\Omega$.
\item $\hbox{dist}_{\rm Euclid}(q, \partial \Omega') > 0$.
\item $q' \in \partial \Omega'$.
\end{itemize}
We exhibit the situation in Figure 1.

\begin{figure}
\centering				     
\includegraphics[height=2in, width=2.75in]{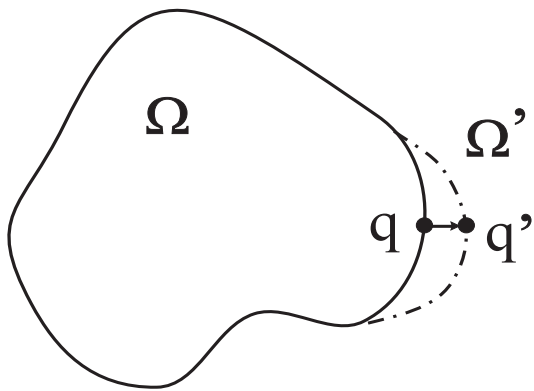}
\caption{The domains $\Omega$ and $\Omega'$.}
\end{figure}

Now let $L_{q'}$ be the Levi polynomial for $\partial \Omega'$ at $q'$. 
Let $\varphi \in C_c^\infty(U)$ be identically equal to 1 in a small
neighborhood of $q'$.

We do not know {\it a priori} that the $\dbar$-problem satisfies uniform estimates on the
domain $\Omega'$.  But we may apply the construction of Beatrous/Range [BEAR] to see that
this is in fact the case (we thank Frank Beatrous and R. Michael Range for helpful remarks regarding this
device).  In detail, suppose $f$ is a $\dbar$ closed form on $\Omega'$. Solve 
$\dbar u = f$ on $\Omega$ with uniform estimates. Let 
$\chi$ be a cutoff function which is 0 on $U$ and identically 1 in the
complement of a slightly larger strongly pseudoconvex neighborhood of $q$.  Let $u_0=\chi u$, extended as zero 
across the perturbed part of the boundary. Let $f_0 = f- \dbar u_0$, which is 
defined and bounded on $\Omega'$ and vanishes in a neighborhood of
$\partial \Omega \setminus U$.  We can therefore solve $\dbar v 
= f_0$ in $\Omega'$, with uniform estimates, by Theorem 1.1 in [BEAR]. 
The solution in $\Omega'$ to the original equation is then 
$u_0+v$.  And that solution is bounded.

Now we use this last result to solve the equation 
$$
\dbar u = (\dbar \varphi) \cdot \frac{1}{L_{q'}} 
$$
on $\Omega'$.  The data on the righthand side 
is $\dbar$-closed and has bounded coefficients.
So there is a bounded solution $u$ by our hypothesis.

Set 
$$
h(z) = \varphi(z) \cdot \frac{1}{L_{q'}(z))} - u , .
$$
Then $h$ is holomorphic and bounded on $\Omega$ and does not
analytically continue past $q'$.  So we may complete
the argument just as in the proofs of Proposition 8.
\endpf
\smallskip \\

\noindent {\bf Corollary of the Proof:}  If $\Omega$ is a smoothly
bounded domain on which uniform estimates for the $\dbar$-equation hold,
and if $\Omega'$ is a domain obtained from $\Omega$ by perturbing
the strongly pseudoconvex points (so that the perturbed points are also
strongly pseudoconvex), then the $\dbar$-problem on $\Omega'$ also
satisfies uniform estimates.
\smallskip \\

We conclude this section by noting that in fact the proof of Theorem 1.1 in [BEAR]
goes through verbatim if ``strongly pseudconvex'' is replaced by ``finite
type in $\CC^2$.   As a result, in view of the discussion above, we have

\begin{proposition} \sl
If $\Omega$ is a smoothly
bounded domain in $\CC^2$ on which uniform estimates for the $\dbar$-equation hold,
and if $\Omega'$ is a domain obtained from $\Omega$ by perturbing
the finite type points (so that the perturbed points are also
finite type), then the $\dbar$-problem on $\Omega'$ also
satisfies uniform estimates.
\end{proposition}

\section{Peak points}

We have seen peak points and peaking functions put to good use in the proof of
Proposition 13.  Now we shall see them in a more general context.

Let $\Omega$ be a domain of type $EL^\infty$.  So $\Omega$ is pseudoconvex, and
there is a strictly larger domain $\widehat{\Omega}$ so that every bounded holomorphic
function on $\Omega$ analytically continues to a bounded holomorphic function $\widehat{f}$ on $\widehat{\Omega}$.
Of course the operator $T: f \mapsto \widehat{f}$ is linear.  It is one-to-one and onto.  It follows
from the closed graph theorem that $T$ continuous.  Now we have a lemma:

\begin{lemma} \sl
The operator $T$ has norm 1.
\end{lemma}
{\bf Proof:}  Of course the norm of $T$ is at least 1.  Suppose that it is actually greater
than 1.  Then there is an $H^\infty$ function $f$ on $\Omega$ so that $f$ has norm 1, and
its extension $\widehat{f}$ has norm greater than 1.  For $k$ a positive integer
consider $g_k = f^k$.  Then the extension of $g_k$ to $\widehat{\Omega}$ is $\widehat{g}_k = (\widehat{f})^k$.
As $k \ra +\infty$, the norm of $\widehat{g}_k$ tends to $+\infty$ while the norm of $g_k$ remains 1.
That is a contradiction.
\endpf
\smallskip \\

\begin{proposition} \sl
Let $\Omega \ss \CC^n$ be a domain and let $q \in \partial \Omega$ be a peak point (see the
proof of Proposition 13).  Let $f_q$ be the peaking function.
Then there cannot be a domain $\widehat{\Omega}$ which properly contains $\Omega$
so that {\bf (i)}  any bounded holomorphic function on $\Omega$ analytically
continues to $\widehat{\Omega}$ and {\bf (ii)}  $q$ lies in the interior of $\widehat{\Omega}$.
\end{proposition}
{\bf Proof:}  Suppose to the contrary that there is such a domain $\widehat{\Omega}$.  Then the
holomorphic function $f_q$ analytically continues to a function $\widehat{f}_q$ on $\widehat{\Omega}$.  Of course
$f_q$ has $H^\infty$ norm 1.  Thus the extended function $\widehat{f}_q$ will also
have norm 1.  But $\widehat{f}_q(q) = 1$.  This contradicts the maximum modulus principle
unless $f_q \equiv 1$.  But that is impossible by the definition of peak function.
\endpf
\smallskip \\

\begin{remark} \rm
In fact one does not need the full force of $q$ being a peak point in order
for this last result to hold.  It is sufficient, for instance, for the
nontangential limit of $f$ at $q$ to be 1, and the values of $f$ at other
points of $\Omega$ have modulus smaller than 1.  

It may also be noted that, by a result of Basener [BAS], the set of peak points
for a domain is contained in the closure of the strongly pseudoconvex points.
This observation is helpful in applying the last proposition.
\end{remark}
	  
In the paper [SIB2], N. Sibony exhibits a smoothly bounded, pseudoconvex
domain on which the equation
$$
\dbar u = f \,,
$$
for $f$ a $\dbar$-closed $(0,1)$ form with bounded coefficents, has
no bounded solution $u$.  This is important information for function
theory, and also for the theory of partial differential equations.

\section{Convexity}

It is natural in the present context to consider the convexity of a domain $\Omega$
with respect to the family of bounded holomorphic functions on $\Omega$.   As usual
we take $\Omega$ to be bounded.

We say that $\Omega$ is {\it convex with respect to the family} ${\cal F} = H^\infty(\Omega)$
if, whenever $K \ss \Omega$ is compact then its {\it hull}
$$
\widehat{K} \equiv \{z \in \Omega: |f(z| \leq \max_{w \in K} |f(w)| \ \hbox{for all} \ f \in {\cal F} \}
$$
is also compact.   When ${\cal F}$ is the family of {\it all} holomorphic functions on $\Omega$ then
$\Omega$ is convex with respect to ${\cal F}$ if and only $\Omega$ is a domain of holomorphy
(see, for instance [KRA1]).   When ${\cal F}$ is the family of $L^2$ holomorphic functions
on $\Omega$ then it is also known (see [PFL]) that $\Omega$ is convex with respect to ${\cal F}$ if
and only if $\Omega$ is a domain of holomorphy.   

We shall now show that $\Omega$ is convex with respect to the family
${\cal F} = H^\infty(\Omega)$ if and only if $\Omega$ is a domain
of type $HL^\infty$.

\begin{proposition}   \sl
If $\Omega$ is convex with respect to $H^\infty(\Omega)$ then $\Omega$ is of type $HL^\infty$.
\end{proposition}
{\bf Proof:}  Assume that $\Omega$ is convex with respect to the family ${\cal F} = H^\infty(\Omega)$.
Now imitate the argument for the implication {\boldmath ${\bf (1)} \Rightarrow {\bf (2)}$}
on pages 154--155 of [KRA1].  Notice that, in that argument, we may replace $h_j$ with
$[1 + h_j]/2$.  Thus each $h_j$ has real part that exceeds $1/4$.  
Notice that 
$$
\left | \prod_{j=1}^\infty (1 - h_j)^j \right |
$$
is bounded (by the standard complex variable theory of infinite products---see [GRK]) by
$$
\exp \left [\sum_j j \cdot 2^{-j} \right ] \, .
$$

Thus we have created a bounded holomorphic function on $\Omega$ that plainly cannot
analytically continue past the boundary.  So $\Omega$ is a domain of type $HL^\infty$.
\endpf 
\smallskip \\

\begin{proposition} \sl
If $\Omega$ is a given domain and if each point of $\partial \Omega$ is essential for $H^\infty(\Omega)$, then
$\Omega$ is convex with respect to $H^\infty(\Omega)$.
\end{proposition}
{\bf Proof:}  Imitate the proof of {\boldmath ${\bf (3)} \Rightarrow {\bf (6)}$} on page 155 of [KRA1].

Now {\boldmath ${\bf (3)} \Rightarrow {\bf (6)}$} on page 155 of [KRA1] shows that if
$f$ is in $H^\infty(\Omega)$ and (for a compact set $K$)
$$
|f(z)| \leq \delta_\Omega(z) \qquad \hbox{for all} \ z \in K
$$
then
$$
|f(z)| \leq \delta_\Omega(z) \qquad \hbox{for all} \ z \in \widehat{K} \, .
$$
Taking $f \equiv \alpha$ for a suitable constant $\alpha < \delta_\Omega(z)$ shows that
$$
\hbox{dist}(K, \partial \Omega) = \hbox{dist}(\widehat{K}, \Omega) \, .
$$
So $\Omega$ is convex with respect to $H^\infty(\Omega)$.
\endpf 
\smallskip \\

\begin{corollary} \sl
If $\Omega$ is of type $HL^\infty$ then $\Omega$ is convex
with respect to the family ${\cal F} = H^\infty(\Omega)$.
\end{corollary}

\section{Concluding Remarks}

It would have been best if we could have given a characterization
of $HL^p$ domains or $EL^p$ domains.  Unfortunately such
a result is beyond our reach at this time.

We hope that the information gathered here will help
to inform the situation and lead, in future work, to
increased understanding of this fascinating problem.
It is clear that there are a spectrum of domains of holomorphy,
and it is in our best interest to understand the elements of this
spectrum.

\newpage

\noindent {\Large \sc References}
\smallskip  \\

\begin{enumerate}

\item[{\bf [BAS]}]  R. Basener, Peak points, barriers, and pseudoconvex boundary points,
{\it Proc.\ Am.\ Math. Soc.} 65(1977), 89-92.

\item[{\bf [BEF]}]  E. Bedford and J.-E. Forn\ae ss, A construction of peak 
functions on weakly pseudoconvex domains, {\it Ann. Math.} 107(1978), 555-568.

\item[{\bf [BER]}]  G. Berg, Bounded holomorphic functions of several variables,
{\it Ark.\ Mat.} 20(1982), 249--270.

\item[{\bf [BERS]}]  L. Bers, {\it Introduction to Several Complex Variables},
New York Univ. Press, New York, 1964.
	
\item[{\bf [BLG]}]  T. Bloom and I. Graham, A geometric characterization of points
of type $m$ on real submanifolds of $\CC^n,$ {\it J. Diff. Geom.} 12(1977), 171-182.

\item[{\bf [CAT1]}]  D. Catlin, Boundary behavior of holomorphic functions on
pseudoconvex domains, {it J. Diff.\ Geometry} 15(1980), 605--625.

\item[{\bf [CAT2]}] D. Catlin, Necessary conditions for subellipticity of
the $\overline{\partial}-$Neumann problem, {\it Ann. Math.} 117(1983),
147-172.

\item[{\bf [CAT3]}] D. Catlin, Subelliptic estimates for the $\dbar$Neumann
problem, {\it Ann. Math.} 126(1987), 131-192.

\item[{\bf [CHS]}]  S.-C. Chen and M.-C. Shaw, {\it Partial Differential Equations
in Several Complex Variables}, American Mathematical Society, Providence, RI, 2001.

\item[{\bf [DAN1]}]  J. P. D'Angelo, Real hypersurfaces, orders of contact,
and applications, {\it Annals of Math.} 115(1982), 615-637.

\item[{\bf [DAN2]}] J. P. D'Angelo, {\it Several Complex Variables and the
Geometry of Real Hypersurfaces}, Studies in Advanced Mathematics. CRC
Press, Boca Raton, FL, 1993.

\item[{\bf [DAR]}]   N. Daras, Existence domains for holomorphic $L^p$ functions,
{\it Publicacions Matem\`{a}tiques} 38(1994), 207--212.

\item[{\bf [FEK]}] C. Fefferman and J. J. Kohn, H\"{o}lder estimates on
domains of complex dimension two and on three dimensional CR-manifolds,
{\it Adv.\ in Math.} 69(1988), 223-303.

\item[{\bf [FOS]}]  J.-E. Forn\ae ss and N. Sibony, On $L^p$ estimates
for $\overline{\partial}$, {\it Several Complex Variables and Complex Geometry}, part 3,
(Santa Cruz, CA, 1989), Proceedings of Symposia in Pure Math., American Mathematical Society,
Providence, RI, 1991, pp.\ 129--163.

\item[{\bf [GRA]}] I. Graham, Boundary behavior of the Carath\'{e}odory and
Kobayashi metrics on strongly pseudoconvex domains in $\CC^n$ with smooth
boundary, {\it Trans. Am. Math. Soc.} 207(1975), 219-240.
    
\item[{\bf [GRK]}]  R. E. Greene and S. G. Krantz, {\it Function Theory of One Complex
Variable}, $2^{\rm nd}$ ed., American Mathematical Society, Providence, RI, 2002.

\item[{\bf [GUR]}]  R. C. Gunning and H. Rossi, {\it Analytic Functions of
Several complex Variables}, Prentice-Hall, Englewood Cliffs, 1965.
						    
\item[{\bf [HAS]}] M. Hakim and N. Sibony, Spectre de
$A(\overline{\Omega})$ pour des domaines born\'{e}s faiblement
pseudoconvexes r\'{e}guliers, {\it Jour.\ Functional Analysis}
37(1980), 127--135.

\item[{\bf [HAP]}]  R. Harvey and J. Polking, Removable singularities of solutions of linear 
partial differential equations, {\it Acta Math.} 125(1970), 39--56. 

\item[{\bf [HEL]}] G. M. Henkin and J. Leiterer, {\it Theory of
Functions on Strictly pseudoconvex Sets with Nonsmooth
Boundary}. With German and Russian summaries. Report MATH
1981, 2. Akademie der Wissenschaften der DDR, Institut für
Mathematik, Berlin, 1981.

\item[{\bf [HEN]}] G. M. Henkin, A uniform estimate for the solution of the
$dbar$-problem in a Weil region, {\it Uspehi Mat.\ Nauk} 26(1971),
211--212.

\item[{\bf [HOR]}] L. H\"{o}rmander, Generators for some rings of analytic
functions, {\it Bull.\ Am.\ Math.\ Soc.} 73(1967), 943--949.

\item[{\bf [KOH]}] J. J. Kohn, Boundary behavior of $\dbar$ on weakly
pseudoconvex manifolds of dimension two, {\em J. Diff. Geom.} 6(1972),
523-542.

\item[{\bf [KON]}] J. J. Kohn and L. Nirenberg, A pseudo-convex domain not
admitting a holomorphic support function, {\it Math. Ann.} 201(1973),
265-268.

\item[{\bf [KRA1]}]  S. G. Krantz, {\it Function Theory of Several complex Variables},
$2^{\rm nd}$ ed., American Mathematical Society, Providence, RI, 2001.

\item[{\bf [KRA2]}] S. G. Krantz, Optimal Lipschitz and $L^{p}$
regularity for the equation $\dbar u = f$ on strongly
pseudo-convex domains, {\it Math. Annalen} 219(1976), 233-260.

\item[{\bf [KRA3]}]  S. G. Krantz, {\it Geometric Function Theory:  Explorations in
Complex Analysis}, Birkh\"{a}user Publishing, Boston, MA, 2006.

\item[{\bf [PFL]}]  P. Pflug, Polynomiale Funktionen auf Steinschen Gebietein Steinschen
Mannigfaltigkeiten, {\it Arch.\ Math.} Basel 28(1977), 169--172.

\newpage

\item[{\bf [SIB1]}]  N. Sibony, Prolongement des fonctions holomorphes born\'{e}es et
m\'{e}trique de Carath\'{e}odory, {\it Inventiones Math.} 29(1975), 205--230.

\item[{\bf [SIB2]}]  N. Sibony, Un exemple de domain pseudoconvexe regulier ou
l'equation $\dbar u = f$ n'admet pas de solution born\'{e}e pour $f$ born\'{e}e,
{\it Invent.\ Math.} 62(1980), 235-242.

\end{enumerate}
\vspace*{.45in}

\noindent \begin{quote}
Department of Mathematics \\
Washington University in St. Louis \\
St.\ Louis, Missouri 63130 \\ 
{\tt sk@math.wustl.edu}
\end{quote}

\end{document}